\documentclass[12pt]{amsart}
\usepackage{mathptmx}
\usepackage{amsmath}
\usepackage{amssymb}
\usepackage{array}
\usepackage{geometry}
\usepackage[bookmarks=true,colorlinks=true, pdfstartview=FitV, linkcolor=black, citecolor=blue, urlcolor=black]{hyperref}
\usepackage{color}
\definecolor{DarkRed}{rgb}{0.55,.00,0.2}
\definecolor{DarkGrey}{rgb}{0.35,.35,0.35}

\theoremstyle{definition}

\theoremstyle{remark}

\numberwithin{equation}{section}

\begin{document}

\title{New inversion, convolution and Titchmarsh's  theorems for the half-Hartley  transform}

\author{S. Yakubovich}
\address{Department of Mathematics, Fac. Sciences of University of Porto,Rua do Campo Alegre,  687; 4169-007 Porto (Portugal)}
\email{ syakubov@fc.up.pt}



\maketitle

\markboth{\rm \centerline{ S.  YAKUBOVICH}}{}
\markright{\rm \centerline{THE HALF-HARTLEY TRANSFORM }}

\begin{abstract}
The generalized Parseval equality for the Mellin transform is employed to prove  the inversion theorem  in $L_2$ with the respective inverse operator related to the Hartley  transform on the nonnegative  half-axis (the half-Hartley transform). Moreover, involving  the convolution method, which is based on the double Mellin-Barnes integrals,  the corresponding
convolution and Titchmarsh's theorems for the half-Hartley transform are established. As an application, we consider solvability conditions for a homogeneous integral equation of the second kind involving the Hartley kernel.   
\bigskip
\end{abstract}

\section{Introduction and auxiliary results}

The familiar reciprocal pair of the Hartley transforms
$$(\mathcal{H} f)(x) = {1\over \sqrt {2\pi} } \int_{-\infty}^\infty [\cos(xt)+ \sin(xt) ] f(t) dt, \ x \in \mathbb{R},\eqno(1.1)$$
$$f(x)= {1\over \sqrt {2\pi} } \int_{-\infty}^\infty [\cos(xt)+ \sin(xt) ] (\mathcal{H} f)(t) dt\eqno(1.2)$$
is well-known \cite{Brace} in connection with various  applications in mathematical physics.  Mapping  and inversion properties of these transforms in $L_2$ as well as their multidimensional analogs were investigated, for instance, in \cite{vu}, \cite{haiw}, \cite{luch}. These operators were treated as the so-called bilateral Watson transform.  Recently, the author found the paper \cite{hart} (see also \cite{half}, \cite{err}), where the attempt to invert the Hartley transform with the integration over $\mathbb{R}_+$
$$(\mathcal{H}_+ f)(x) = \sqrt { {2\over \pi} } \int_{0}^\infty [\cos(xt)+ \sin(xt) ] f(t) dt, \ x \in \mathbb{R}_+,\eqno(1.3)$$
was  undertaken. However, the inversion formula obtained by the authors is depending on the Fourier transform of the image and, indeed, needs to be improved. Here we will achieve this main goal, proving the inversion theorem for transformation (1.3) in $L_2(\mathbb{R}_+)$. Moreover, we will construct and study properties of the convolution operator, related to the half-Hartley transform by general convolution method developed  by the author in 1990,  and  which is based on the double Mellin-Barnes integrals  \cite{con}, \cite{hai}, \cite{luch}. Namely, we will prove the convolution theorem and Titchmarsh's theorem about the absence of divisors in the convolution product. Finally, we apply the half-Hartley transform (1.3) to find solvability conditions and the form of solutions for  a homogeneous integral equation of the second kind.

We note in this section that our natural approach will involve the $L_2$-theory of the Mellin transform \cite{tit}
$$(\mathcal {M} f)(s)= f^*(s)= \int_0^\infty f(t) t^{s-1}dt, \ s \in
\sigma  =\{ s \in \mathbb{C}, s={1\over 2}  +i\tau\},\eqno(1.4)$$
where the integral  is convergent in the mean square sense with respect to the norm in $L_2(\sigma)$. Reciprocally,  the inversion formula takes place
$$f(x)= {1\over 2\pi i}\int_{\sigma} f^*(s) x^{-s} ds,\ x >0\eqno(1.5)$$
with the convergence of the integral in the mean square sense with respect to the norm in $L_2(\mathbb{R}_+)$.  Furthermore, for any $f_1, f_2  \in L_2(\mathbb{R}_+)$ the generalized Parseval identity holds
$$\int_0^\infty f_1\left(xt\right) f_2(t) dt  = {1\over 2\pi i}\int_{\sigma}  f_1^*(s)f_2^*(1-s) x^{-s}
ds, \ x >0\eqno(1.6)$$
with  Parseval's equality of squares of $L_2$- norms
$$\int_0^\infty |f(x)|^2 dx = {1\over 2\pi}   \int_{-\infty}^{\infty} \left|f^*\left({1\over 2}  + i\tau\right)\right|^2 d\tau.\eqno(1.7)$$

\section{Inversion theorem}

We begin with the following  inversion theorem for the half-Hartley  transform (1.3).  Precisely, it has

{\bf Theorem 1}. {\it  The half-Hartley transform $(1.3)$ extends to a bounded invertible map
$\mathcal{H}_+:  L_2(\mathbb{R}_+) \to L_2(\mathbb{R}_+)$ and for almost all $x \in \mathbb{R}_+$ the following reciprocal formulas hold
$$(\mathcal{H}_+ f)(x) = \sqrt { {2\over \pi} } {d\over dx} \int_{0}^\infty [1+\sin(xt) - \cos (xt) ] {f(t)\over t} dt, \ x \in \mathbb{R}_+,\eqno(2.1)$$
$$f(x)=   \sqrt { {2\over \pi} }  \int_0^\infty \left[ \sin( xt)  \ S(xt)+ \cos (xt) \  C(xt)\right] (\mathcal{H}_+ f)(t) dt,\eqno(2.2)$$
where $S(x),\ C(x)$ are Fresnel sin- and cosine- integrals, respectively,
$$S(x)=  \sqrt { {2\over \pi} } \int_0^{\sqrt x}  \sin(t^2) dt, \quad   C(x)=  \sqrt { {2\over \pi} } \int_0^{\sqrt x}  \cos (t^2) dt$$
and  integral $(2.2)$ converges with respect to the norm in $L_2(\mathbb{R}_+)$.  Finally,   the norm inequalities take place}
$$\sqrt 2 \left|\left| f \right|\right|_{ L_2(\mathbb{R}_+)}  \le \left|\left|  \mathcal{H}_+ f \right|\right|_{ L_2(\mathbb{R}_+)} \le
 2   \left|\left| f \right|\right|_{ L_2(\mathbb{R}_+)}.\eqno(2.3)$$

\begin{proof}  Let $f$ belong to the space $C^{(2)}_c(\mathbb{R}_+)$  of  continuously differentiable functions
of compact support, which is dense in $ L_2(\mathbb{R}_+)$.  Then integrating by parts  in (1.4), we find that $s^2 f^*(s)$ is bounded on $\sigma$ and therefore $ f^*(s) \in L_2(\sigma) \cap L_1(\sigma) $.   Hence minding the known formulas \cite{tit}
$$\int_0^\infty {\sin t \over t} t^{s-1} dt = \frac{\Gamma(s)}{1-s} \cos\left({\pi s\over 2}\right), \ s \in \sigma,\eqno(2.4)$$
$$\int_0^\infty {1-\cos t \over t} t^{s-1} dt = \frac{\Gamma(s)}{1-s} \sin \left({\pi s\over 2}\right), \ s \in \sigma,\eqno(2.5)$$
we call the generalized Parseval equality (1.6) to derive for all $x >0$
$$\sqrt { {2\over \pi} } \int_{0}^\infty [1+\sin(xt) - \cos (xt) ] {f(t)\over t} dt$$$$ =  \sqrt { {2\over \pi} } \ {1\over 2\pi i}\int_{\sigma}  \Gamma(s)\left[ \sin \left({\pi s\over 2}\right) + \cos\left({\pi s\over 2}\right)\right] f^*(1-s)  \frac{x^{1 -s}}{1-s}  ds.\eqno(2.6)$$
It is easily seen the possibility to differentiate through with respect to $x$ in equality (2.6).  Thus combining with (1.3), we derive (2.1) together with the equality
$$ (\mathcal{H}_+ f)(x) = \sqrt { {2\over \pi} } \ {1\over 2\pi i}\int_{\sigma}  \Gamma(s)\left[ \sin \left({\pi s\over 2}\right) + \cos\left({\pi s\over 2}\right)\right] f^*(1-s)  x^{-s} ds, \eqno(2.7)$$
which is valid for any $f \in C^{(2)}_c(\mathbb{R}_+)$.  Furthermore, from (1.7) one immediately obtains the norm estimates
$$\left|\left| \mathcal{H}_+ f \right|\right|_{ L_2(\mathbb{R}_+)} =   {2\over \sqrt \pi }  \left( {1\over 2\pi}
 \int_{-\infty}^{\infty} \cosh^2 \left({\pi \tau\over 2}\right)\left| \Gamma\left({1\over 2} + i\tau\right) f^*\left({1\over 2} + i\tau\right)\right|^2 d\tau\right)^{1/2} $$
 $$= 2 \left( {1\over 2\pi}  \int_{-\infty}^{\infty} \frac{\cosh^2 \left(\pi \tau/2 \right)}{\cosh(\pi\tau)} \left|f^*\left({1\over 2} + i\tau\right)\right|^2 d\tau\right)^{1/2} $$
$$=  2 \left( {1\over 2\pi}  \int_{-\infty}^{\infty} \frac{\cosh^2 \left(\pi \tau/2 \right)}{2 \cosh^2(\pi\tau/2)-1} \left|f^*\left({1\over 2} + i\tau\right)\right|^2 d\tau\right)^{1/2} \le 2 \left|\left| f \right|\right|_{ L_2(\mathbb{R}_+)}$$
and plainly
$$ \left|\left| \mathcal{H}_+ f \right|\right|_{ L_2(\mathbb{R}_+)} \ge \sqrt 2 \left|\left| f \right|\right|_{ L_2(\mathbb{R}_+)}.$$
Thus  we proved (2.3) for  any $f \in C^{(2)}_c(\mathbb{R}_+)$.   Further, since $C^{(2)}_c(\mathbb{R}_+)$ is dense in $ L_2(\mathbb{R}_+)$, there is a unique extension of $\mathcal{H}_+$ as an invertible continuous map $\mathcal{H}_+:  L_2(\mathbb{R}_+) \to L_2(\mathbb{R}_+)$.
Now, let $f \in L_2(\mathbb{R}_+)$. There is a sequence $\{ f_n\},\  f_n \in  C^{(2)}_c(\mathbb{R}_+)$
such that $\left|\left| f_n- f  \right|\right|_{ L_2(\mathbb{R}_+)}  \to 0, \ n \to \infty$.  Denoting by
$$h_n(x)= \sqrt { {2\over \pi} } \ {1\over 2\pi i}\int_{\sigma}  \Gamma(s)\left[ \sin \left({\pi s\over 2}\right) + \cos\left({\pi s\over 2}\right)\right] f_n^*(1-s)  x^{ -s} ds, \eqno(2.8)$$
we observe by virtue of (2.3)  that  $\{h_n\}$ is a Cauchy sequence and it has a limit in  $L_2(\mathbb{R}_+)$, which we will call $h$.  Hence,  integrating through in (2.8), we have
$$\int_0^x h_n(y) dy =  \sqrt { {2\over \pi} }\int_0^x \left( \ {1\over 2\pi i}\int_{\sigma}  \Gamma(s)\left[ \sin \left({\pi s\over 2}\right) + \cos\left({\pi s\over 2}\right)\right] f_n^*(1-s)  y^{ -s} ds\right) dy.\eqno(2.9)$$
In the meantime,   by the Schwarz inequality
$$\int_0^x  [ h_n(y)- h(y)] dy \le  \sqrt x  \   \left|\left| h_n- h \right|\right|_{ L_2(\mathbb{R}_+)} \to 0,\ n \to \infty$$
 and in the right-hand side of (2.9) one can change the order of integration by Fubini's theorem. Then passing  to the limit when $n \to \infty$ under integral signs in the obtained equality due to the Lebesgue dominated convergence theorem, we   find
$$\int_0^x h (y) dy =  \sqrt { {2\over \pi} } \ {1\over 2\pi i}\int_{\sigma}  \Gamma(s)\left[ \sin \left({\pi s\over 2}\right) + \cos\left({\pi s\over 2}\right)\right] f^*(1-s)  \frac{x^{1 -s}}{1-s}  ds.\eqno(2.10)$$
Differentiating by $x$ in (2.10),  we come out with the equality for almost all $x>0$
$$ h (x) \equiv (\mathcal{H}_+f)(x)=   \sqrt { {2\over \pi} } \ {1\over 2\pi i}{d\over dx} \int_{\sigma}  \Gamma(s)\left[ \sin \left({\pi s\over 2}\right) + \cos\left({\pi s\over 2}\right)\right] f^*(1-s)  \frac{x^{1 -s}}{1-s}  ds,\eqno(2.11)$$
which coincides with (2.7) for any $f \in C^{(2)}_c(\mathbb{R}_+)$. Consequently,  appealing to (2.4), (2.5),  (1.6) and (2.6),  we complete the proof of  representation (2.1).

 Finally, we establish the inversion formula (2.2).    To do this, we denote by $h^*(s)$ the Mellin transform of $h(t)$ (1.4) in $L_2$ and write reciprocally to (2.11) for almost all $x >0$
$$ f (x) =   \sqrt { {\pi\over 2} } \ {1\over 2\pi i}{d\over dx} \int_{\sigma}  \Gamma^{-1} (1-s)\left[ \sin \left({\pi s\over 2}\right) + \cos\left({\pi s\over 2}\right)\right]^{-1} h^*(1-s)  \frac{x^{1 -s}}{1-s}  ds.$$
Meanwhile with the supplement formula for gamma-functions and elementary trigonometric manipulations it becomes
$$ f (x) =    {1\over\sqrt { 2\pi} }  \ {1\over 2\pi i}{d\over dx} \int_{\sigma}  \Gamma (s)\frac{\sin(\pi s) \left[ \sin \left(\pi s/2\right) + \cos\left(\pi s/ 2\right)\right]}{1+ \sin(\pi s)}  h^*(1-s)  \frac{x^{1 -s}}{1-s}  ds$$
$$= {1\over 2}   (\mathcal{H}_+h)(x) - {1\over\sqrt { 2\pi} }  \ {1\over 2\pi i}{d\over dx} \int_{\sigma}  \Gamma (s)\frac{ \sin \left(\pi s/2\right) + \cos\left(\pi s/ 2\right)}{1+ \sin(\pi s)}  h^*(1-s)  \frac{x^{1 -s}}{1-s}  ds$$
$$= {1\over 2}   (\mathcal{H}_+h)(x) - {1\over\sqrt { 2\pi} }  \ {1\over 2\pi i} \int_{\sigma}  \frac{\Gamma (s) }{\sin \left(\pi s/2\right) + \cos\left(\pi s/ 2\right)}  h^*(1-s) x^{ -s} ds,$$
$$= {1\over 2}   (\mathcal{H}_+h)(x) - {1\over2 \pi \sqrt {\pi} }  \ {1\over 2\pi i} \int_{\sigma}  \Gamma (s)
\Gamma\left({s\over 2}+ {1\over 4} \right) \Gamma\left({3\over 4} - {s\over 2} \right) h^*(1-s) x^{ -s} ds,$$
where the differentiation under the integral sign is allowed via the absolute and uniform convergence.  But the generalized Parseval identity (1.6)  yields
$$f(x)=  {1\over 2}   (\mathcal{H}_+h)(x) - {1\over2 \pi \sqrt {\pi} } \int_0^\infty k(xt) h(t) dt,\eqno(2.12)$$
where
$$k(x)= {1\over 2\pi i} \int_{\sigma}  \Gamma (s) \Gamma\left({s\over 2}+ {1\over 4} \right) \Gamma\left({3\over 4} - {s\over 2} \right) x^{ -s} ds,\ x >0.$$
Hence employing again (1.6) and using relations  (8.4.2.5), (8.4.3.1)  in \cite{prud}, Vol. 3, the latter integral can be written in the form
$$k(x)= 2 \int_0^\infty \frac{ e^{-x/t}  }{ 1+ t^2}\  {dt\over \sqrt t} =
2 \int_0^\infty \frac{ e^{-x t} \sqrt t  }{ 1+ t^2} dt $$
and it is calculated via relation (2.3.7.10) in \cite{prud}, Vol. 1, namely,
$$k(x)= \pi \sqrt 2 [\sin x + \cos x]-  2^{3/2} \pi \left[ \sin x \ S(x)+ \cos x\  C(x)\right],$$
where  $S(x),\ C(x)$ are Fresnel sin- and cosine- integrals (see above). Hence, substituting the value of $k(x)$ in (2.12) and making use (2.1), we come out with the inversion formula (2.2) and  and complete the proof of Theorem 1.
\end{proof}

{\bf Remark 1}.  Taking into account the value of the integral (2.5.5.1) in \cite{prud}, Vol. 1 and relation (7.14.2.75) in \cite{prud}, Vol. 3, we have
$$  \sqrt { {2\over \pi} }  \left[ \sin( x)  \ S(x)+ \cos (x) \  C(x)\right] =
 {1\over \pi}\int_0^x { \cos t \over \sqrt {x- t} } dt $$$$=  {1 \over 2\sqrt \pi}
\left[ e^{i(x- \pi/4)} \hbox{erf} \left( e^{i\pi/4} \sqrt x \right)+   e^{- i(x+ \pi/4)} \hbox{erfi} \left( e^{i\pi/4} \sqrt x \right)\right],$$
where $\hbox{erf} (z)$ and $\hbox{erfi} (z)$ are the error function and the error function of imaginary argument, respectively,
$$\hbox{erf} (z)= {2\over \sqrt \pi} \int_0^z e^{-t^2} dt, \quad   \hbox{erfi} (z)= {2\over \sqrt \pi} \int_0^z e^{t^2} dt.$$
Therefore,   inversion formula (2.2)  can be written as
$$f(x)= {1 \over 2\sqrt \pi}  \int_0^\infty
\left[ e^{i(xt - \pi/4)} \hbox{erf} \left( e^{i\pi/4} \sqrt {xt} \right)+   e^{- i(xt+ \pi/4)} \hbox{erfi} \left( e^{i\pi/4} \sqrt {xt} \right)\right]
(\mathcal{H}_+ f)(t) dt.$$

\section{Convolution operator for the half-Hartley transform}

In this section we will construct and study mapping properties of  the convolution, related to the transformation (1.3).
 Following the general convolution method developed    for integral transforms of the Mellin convolution type
  (cf. \cite{con}, \cite{hai},  \cite{luch})  we have

{\bf Definition 1}. Let   $f, g$ be functions from $\mathbb{R}_+$ into $\mathbb{C}$ and $f^*,\  g^*$ be their Mellin transforms $(1.4)$.
 Then the function $f*g$ being defined on $\mathbb{R}_+$ by the double Mellin-Barnes integral
$$(f*g)(x)= \frac{1}{ ( 2\pi i)^2} \int_{\sigma}  \int_{\sigma} \frac{\Gamma(s)\Gamma(w)}{\Gamma (s+w-1/2)}
 \frac {\sin(\pi(s+w)/2) + \cos(\pi(s-w)/2)} {\sin(\pi(s+w)/2) } $$$$\times   f^*(1-s)g^*(1-w) x^{s+w- 3/2}dsdw\eqno(3.1)$$
is called the convolution of $f$ and $g$  related to the half-Hartley transform $(1.3)$  (provided that it exists).

The  convolution theorem for the half-Hartley  transform can be stated as

 {\bf Theorem 2}.  {\it  Let  $f^*,\ g^*$ be the Mellin transforms of $f, g$, respectively, satisfying
  conditions $s f^*(s),\\  s g^*(s)  \in  L_2(\sigma)$.   Then the convolution $(3.1)$ $f*g$ exists and belongs
   to $L_2(\mathbb{R}_+)$ with the norm estimate
$$||f*g||_{ L_2(\mathbb{R}_+)} \le  4\sqrt {{2\over \pi}} \left( \int_{-\infty}^\infty  \left |(1/2+i\theta) g^*(1/2+i\theta)\right|^2 d\theta \right)^{1/2}
 $$$$\times  \left(\int_{-\infty}^\infty \left|(1/2+i\tau)  f^*(1/2+ i \tau)\right|^2 d\tau \right)^{1/2} d\theta.\eqno(3.2)$$
 Moreover, its  Mellin transform  $(\mathcal { M} \  ( f*g) ) (1- s) $  is equal to
 $$(\mathcal { M}  (f*g) ) (1-s) =   \frac{\sqrt 2 }{2\pi  i \  \Gamma(s)   \sin(\pi( s+1/2)/2)}  \int_{\sigma}
  \Gamma(s-w+1/2)   \cos(\pi( s-w)/2) \Gamma(w) $$$$\times   \sin(\pi( w+ 1/2)/2)
   f^*(1/2-s+w )g^*(1-w)dw,\   s \in \sigma.\eqno(3.3)$$
Besides,  the factorization equality  holds
$$(\mathcal{H}_+   (f*g) )(x) =   \sqrt {{x\pi \over 2}} \ (\mathcal{H}_+ f)(x) (\mathcal{H}_+ g)(x), \quad x >0 \eqno(3.4)$$
as well as the generalized Parseval type identity
$$(f*g)(x)=   \int_0^\infty \left[ \sin( xt)  \ S(xt)+ \cos (xt) \
C(xt)\right]\sqrt t\  (\mathcal{H}_+ f)(t)(\mathcal{H}_+ g)(t)
dt,\eqno(3.5)$$
where integral $(3.5)$ converges in the $L_2$- sense.}

\begin{proof}   In fact, for $s= 1/2 + i\tau, w= 1/2+ i\theta, \ (\tau, \theta ) \in \mathbb{R}$ we obtain
$$\left| \frac{\Gamma(s)\Gamma(w)}{\Gamma (s+w-1/2)} \frac {\sin(\pi(s+w)/2) + \cos(\pi(s-w)/2)} {\sin(\pi(s+w)/2)
}\right| = 2\left|\frac{\Gamma(1/2+i\tau)\Gamma(1/2+i\theta)}{\Gamma
(1/2+ i(\tau+\theta))}\right| $$$$\times \frac
{\cosh(\pi\tau/2)\cosh(\pi\theta/2)} {\cosh(\pi(\tau+\theta)/2)}= 2
\sqrt \pi \frac {\cosh^{1/2}(\pi(\tau +
\theta))\cosh(\pi\tau/2)\cosh(\pi\theta/2)}
{\cosh^{1/2}(\pi\tau)\cosh^{1/2}(\pi\theta)\cosh(\pi(\tau+\theta)/2)}$$$$
\le 2 \sqrt {2\pi} \frac {\cosh(\pi\tau/2)\cosh(\pi\theta/2)}
{\cosh^{1/2}(\pi\tau)\cosh^{1/2}(\pi\theta)}\le 4\sqrt {2\pi}.\eqno(3.6)$$
Meanwhile, plainly via conditions $s f^*(s),  s g^*(s)  \in  L_2(\sigma)$ we have
 $f^*(s),\  g^*(w) \in L_1(\sigma)$. Therefore by Fubini's theorem  the double integral (3.1) is equal to the corresponding iterated integrals.
 Then  making in (3.1) the simple substitution $z= s+w-1/2$, using elementary trigonometric formulas,  (1.5), (1.7) and the above estimate,
 we easily come out with (3.3) and  the estimate (3.2). Namely,
 it  has
$$||f*g||_{ L_2(\mathbb{R}_+)}  = \left({1\over 2\pi} \int_{-\infty}^\infty \left|(\mathcal { M}  (f*g) ) (1/2+ i\tau)\right|^2 d\tau\right)^{1/2}
 $$$$ \le {2\over \pi} \int_{-\infty}^\infty  \left |g^*(1/2+i\theta)\right|
 \left(\int_{-\infty}^\infty \left| f^*(1/2+ i (\tau-\theta)\right|^2 d\tau \right)^{1/2} d\theta
  \le {4\over \pi} \left( \int_{-\infty}^\infty  {d\theta \over  \theta^2+ 1/4} \right)^{1/2}$$$$
  \times \left( \int_{-\infty}^\infty  \left |(1/2+i\theta) g^*(1/2+i\theta)\right|^2 d\theta \right)^{1/2}
   \left(\int_{-\infty}^\infty \left|(1/2+i\tau)  f^*(1/2+ i \tau)\right|^2 d\tau \right)^{1/2} d\theta $$
$$= 4\sqrt {{2\over \pi}} \left( \int_{-\infty}^\infty  \left |(1/2+i\theta) g^*(1/2+i\theta)\right|^2 d\theta \right)^{1/2}
  \left(\int_{-\infty}^\infty \left|(1/2+i\tau)  f^*(1/2+ i \tau)\right|^2 d\tau \right)^{1/2} d\theta,$$
where the Schwarz and generalized Minkowski inequalities are
employed. The factorization equality (3.4) comes immediately from
(3.3) and (1.6) because
$$(\mathcal{H}_+  (f*g) )(x) =  {2\over \sqrt \pi}\ \frac{1  }{2\pi  i} \int_{\sigma}  \
 (\mathcal { M}  (f*g) ) (1-z)  \Gamma(z)   \sin(\pi( z+ 1/2)/2) x^{-z} dz$$$$=
\sqrt{{2\over \pi}}\  \frac{ 2}{ ( 2\pi i)^2} \int_{\sigma}
\int_{\sigma} \Gamma(z-w+1/2)   \cos(\pi( z-w)/2)
\Gamma(w)\sin(\pi(w+ 1/2)/2)$$$$\times  f^*(1/2-z+w )g^*(1-w) x^{-z}
dw dz$$$$= \sqrt{{2\over \pi}}\  \frac{ 2\sqrt x}{ ( 2\pi i)^2}
\int_{\sigma} \int_{\sigma} \Gamma(s)   \cos(\pi( (s-1/2)/2)
\Gamma(w)\sin(\pi(w+ 1/2)/2)$$$$\times  f^*(1-s)g^*(1-w) x^{-s-w} dw
dz = \sqrt {{x\pi \over 2}} \ (\mathcal{H}_+ f)(x) (\mathcal{H}_+
g)(x).$$
Finally, the generalized Parseval identity (3.5) is a direct
consequence of the inversion formula (2.2).

\end{proof}

This section ends with an analog of the Titchmarsh theorem about the
absence of divisors of zero in  convolution (3.1). We have

{\bf Theorem 3}. {\it Let $f^*,\ g^*$ satisfy conditions  $e^{\pi  |s|}  f^*(s),\  e^{\pi  |s|}  g^*(s)  \in
L_1(\sigma)$.   Then if  $(f*g)(x)= 0,  \ x >0$, then either $f(x)=0$ or $g(x) =0$ on $\mathbb{R}_+$.}

\begin{proof}  In fact, the integral
$$F(z)= \frac{1}{ ( 2\pi i)^2} \int_{\sigma}  \int_{\sigma} \frac{\Gamma(s)\Gamma(w)}{\Gamma (s+w-1/2)}
 \frac {\sin(\pi(s+w)/2) + \cos(\pi(s-w)/2)} {\sin(\pi(s+w)/2) } $$$$\times   f^*(1-s)g^*(1-w) z^{s+w- 3/2}dsdw$$
represents an analytic function in the domain $D = \{ z \in \mathbb{C}:  \  |\arg z| < \pi\}$, since under condition of the theorem it converges uniformly for any $z \in \mathbb{C}:  \  |z| \ge a >0, \  |\arg z| < \pi$. Precisely, we have $(s= 1/2+i\tau,\  w= 1/2 +i\theta,\ z^{s+w-3/2}= |z|^{-1/2} e^{(\tau+\theta)\arg z }$) via (3.6)
$$\int_{\sigma}  \int_{\sigma} \left|  \frac{\Gamma(s)\Gamma(w)}{\Gamma (s+w-1/2)}
 \frac {\sin(\pi(s+w)/2) + \cos(\pi(s-w)/2)} {\sin(\pi(s+w)/2) } \right.$$$$\times  \left. f^*(1-s)g^*(1-w) z^{s+w- 3/2}dsdw \right| $$$$\le 4\  \sqrt{{2\pi\over a}} \int_{-\infty}^\infty  \int_{-\infty}^\infty e^{\pi [  |\tau| + |\theta| ]}
  \left| f^*(1/2+i\tau) g^*(1/2 + i\theta)\right|  d\tau d\theta < \infty.$$
Moreover,  (3.1) yields that $F(x)= (f*g)(x)$. Thus by virtue the uniqueness theorem for analytic functions $F(z)= (f*g)(z),
 z \in D$. Moreover, equality (3.4) holds for  $z \in D$, where the main branch of the square root is chosen.
 Hence calling (2.7), we deduce
 $$\left| (\mathcal{H}_+ f)(z)\right| \le  {1\over \pi\sqrt{2\pi}}  \int_{\sigma} \left| \Gamma(s)\left[ \sin \left({\pi s\over 2}\right) + \cos\left({\pi s\over 2}\right)\right] f^*(1-s)  z^{-s} ds\right| $$$$ \le  {1\over \pi\sqrt a } \int_{-\infty}^{\infty} \frac{\cosh(\pi\tau/2)}{\cosh^{1/2}(\pi\tau)}  \left| f^*(1/2+i\tau)\right|   e^{\pi |\tau|} d\tau \le  {1\over \pi} \sqrt {{2\over a} } \int_{-\infty}^{\infty}   \left| f^*(1/2+i\tau)\right|   e^{\pi |\tau|} d\tau < \infty,$$
which means that $ (\mathcal{H}_+ f)(z)$ is analytic in $ D$.   Therefore, if $(f*g)(x)=0, \ x >0,$ then via the uniqueness theorem $(f*g)(z)\equiv 0,\ z \in D$  and (3.4) yields
$$ (\mathcal{H}_+ f)(z) (\mathcal{H}_+ g)(z) =0,\quad z \in D.$$
Since the left-hand side of the latter equality is the product of analytic functions in  $D$,
it means that either  $ (\mathcal{H}_+ f)(z)  \equiv 0$, or $ (\mathcal{H}_+ g)(z) \equiv 0$  in  $D$.
 Finally, we observe that under conditions of the theorem $f, g \in  L_2(\mathbb{R}_+) $ and from operational properties of the inverse Mellin transform (1.5)  it follows that $f, g$ are infinite times differentiable functions.   Thus  appealing to Theorem 1 and inversion formula (2.2) we find   that either $f=0$ or $g=0$ on $\mathbb{R}_+$.
\end{proof}

\section{A homogeneous integral equation involving the  half-Hartley transform}

Here  we will establish solvability  conditions in  $L_2(\mathbb{R}_+)$,  concerning the  following  integral equation of the second kind 
$$  \sqrt { {2\over \pi} } \int_{0}^\infty [\cos(xt)+ \sin(xt) ] f(t) dt=\lambda f(x), \ x \in \mathbb{R}_+,\ \lambda \in \mathbb{C}.\eqno(4.1)$$
The main result of the section is 

{\bf Theorem 4}.  {\it Let $|\lambda | < \sqrt 2$.  In  order to an arbitrary function $f \in L_2(\mathbb{R}_+)$ be a solution of integral equation $(4.1)$  it is necessary to have the form of the integral in the mean square sense  
$$f(x)= {1\over 2\pi i} \int_\sigma \left[ \lambda +  \sqrt{{\pi\over 2} }  \frac{1} {\Gamma(1-s)} \left[\sec\left({\pi s\over 2}\right)+ \csc\left({\pi s\over 2}\right)\right]\right]  \varphi (s) x^{-s} ds,  \  x> 0\eqno(4.2)$$
in terms of some function $\varphi(s)$, satisfying  condition $\varphi(s)= \varphi (1-s),\ s \in \sigma$, i.e.  $\varphi (1/2+i\tau)$ is even with respect to $\tau \in \mathbb{R}$.   Besides,    its Mellin's   transform 
$$f^*(s)=  \left[ \lambda +  \sqrt{{\pi\over 2} }  \frac{1}{\Gamma(1- s)}\left[\sec\left({\pi s\over 2}\right)+ \csc\left({\pi s\over 2}\right) \right] \right] \varphi (s) , \ s \in \sigma \eqno(4.3)$$
as well as  $\varphi$ belong to $L_2(\sigma)$ and the following $L_2$- norm estimates hold
$$ (2+ |\lambda|)^{-1}   ||f^* ||_{L_2(\sigma)}  \le ||\varphi ||_{L_2(\sigma)} \le \left(\sqrt 2 - |\lambda|\right)^{-1}  ||f^* ||_{L_2(\sigma)}.\eqno(4.4)$$ 
The condition $\varphi(s)= \varphi (1-s),\ \varphi \in L_2(\sigma)$ and the form of solutions $(4.2)$ are also sufficient for any $\psi$ as the inverse Mellin transform of $\varphi$, satisfying the integral equation
$$(\lambda^2-2) \psi (x)-   {2\over \pi} \int_0^\infty {\psi(t)\over x+t} dt= 0,  \  x \in  \mathbb{R}_+,\eqno(4.5)$$
where integral $(4.5)$ converges absolutely.}

\begin{proof}  Let $f$ be a solution of equation (4.1).  Then via formulas (2.1) and (2.6) we derive the equality
$$\lambda \int_0^x f(y)dy =   \sqrt { {2\over \pi} } \ {1\over 2\pi i}\int_{\sigma}  \Gamma(s)\left[ \sin \left({\pi s\over 2}\right) + \cos\left({\pi s\over 2}\right)\right] f^*(1-s)  \frac{x^{1 -s}}{1-s}  ds.\eqno(4.6)$$
In the meantime,  its left-hand side is equal to (see (1.6))
$$\lambda \int_0^x f(y)dy = {\lambda \over 2\pi i} \int_{\sigma}   f^*(s) { x^{1-s}\over 1- s } ds.\eqno(4.7)$$
Hence, comparing  right-hand sides  of (4.6), (4.7)  and since both integrand are from $L_1(\sigma)$ they are equal by virtue of Th. 32 in \cite{tit}.   Hence
$$\lambda f^*(s)= \sqrt { {2\over \pi} } \  \Gamma(s)\left[ \sin \left({\pi s\over 2}\right) + \cos\left({\pi s\over 2}\right)\right] f^*(1-s),  \  s \in \sigma.\eqno(4.8)$$ 
But $1-s \in \sigma$.   Therefore,  changing $s$ on $1-s$ in (4.8),  it  becomes
$$\lambda f^*(1- s)= \sqrt { {2\over \pi} } \  \Gamma(1- s)\left[ \sin \left({\pi s\over 2}\right) + \cos\left({\pi s\over 2}\right)\right] 
f^*(s).\eqno(4.9)$$ 
Subtracting (4.9) from (4.8) and then using the supplement and duplication formulas for gamma - functions,  we obtain 
$$\lambda[ f^*(s)-  f^*(1-s)]  = \sqrt  {2\pi}   \left[ \frac{1}{\Gamma(s/2) \Gamma(1- s/2) }  + \frac{1}{\Gamma((1-s)/2) \Gamma((1+ s)/2) }\right] $$$$\times \left[  \Gamma(s) f^*(1-s) -   \Gamma(1-s) f^*(s)\right]=  2^{s-1/2}  f^*(1-s)  \left[ \frac{ \Gamma((1+s) /2)}{ \Gamma(1- s/2) }  + \frac{ \Gamma(s/2)}{\Gamma((1-s)/2) }\right] $$$$- 2^{1/2-s}  f^*(s) \left[ \frac{ \Gamma((1-s) /2)}{ \Gamma(s/2) }  + \frac{ \Gamma(1- s/2)}{\Gamma((1+s)/2) }\right].$$
Hence
$$f^*(s) \left[ \lambda +  2^{1/2-s}  \left[ \frac{ \Gamma((1-s) /2)}{ \Gamma(s/2) }  + \frac{ \Gamma(1- s/2)}{\Gamma((1+s)/2) }\right]\right]
$$$$= f^*(1-s) \left[ \lambda + 2^{s-1/2}  \left[ \frac{ \Gamma((1+s) /2)}{ \Gamma(1- s/2) }  + \frac{ \Gamma(s/2)}{\Gamma((1-s)/2) }\right] \right] ,\  s \in \sigma,\eqno(4.10)$$
or,  since
$$\sqrt 2  \le \sqrt{{\pi\over 2} } \left| \frac{1} { \Gamma(s)} \left[\sec\left({\pi s\over 2}\right)+ \csc\left({\pi s\over 2}\right) \right]\right| =   \frac{2 \cosh(\pi\tau/2)}{\cosh^{1/2} (\pi\tau)} \le 2 , \ s= 1/2 +i\tau,\ \tau \in \mathbb{R},\eqno(4.11)$$
(4.10) under condition $|\lambda| < \sqrt 2$ yields
$$f^*(s) \left[ \lambda +  \sqrt{{\pi\over 2} }  \frac{1} {\Gamma(1-s)} \left[\sec\left({\pi s\over 2}\right)+ \csc\left({\pi s\over 2}\right)\right]\right]^{-1} $$$$= f^*(1-s) \left[ \lambda + \sqrt{{\pi\over 2} }  \frac{1} {\Gamma(s)} \left[\sec\left({\pi s\over 2}\right)+ \csc\left({\pi s\over 2}\right)\right] \right] ^{-1} = \varphi (s),\  s \in \sigma.\eqno(4.12)$$
Thus  we find  that $\varphi (s)= \varphi (1-s)$, i.e is even with respect to $\tau \in \mathbb{R}$, where $s= 1/2 + i\tau$.  
So, we established (4.3)  and reciprocally in $L_2$ the representation (4.2).  Meanwhile estimates (4.11) yield (4.4).

Now, let us assume  for some $\varphi \in L_2(\sigma)$ in (4.3) the condition $\varphi(s)=\varphi(1-s),\ s \in \sigma$. Then  substituting the value of $f^*(s)$ in (4.3) into equation (4.8), we obtain

$$\lambda  \left[ \lambda +  \sqrt{{\pi\over 2} }  \frac{1}{\Gamma(1- s)}\left[\sec\left({\pi s\over 2}\right)+ \csc\left({\pi s\over 2}\right) \right] \right] \varphi (s) =  \sqrt { {2\over \pi} } \  \Gamma(s)\left[ \sin \left({\pi s\over 2}\right) + \cos\left({\pi s\over 2}\right)\right] $$$$\times  \left[ \lambda +  \sqrt{{\pi\over 2} }  \frac{1}{\Gamma( s)}\left[\sec\left({\pi s\over 2}\right)+ \csc\left({\pi s\over 2}\right) \right] \right] \varphi (s), \ s \in \sigma,$$
or after simple calculations it drives  to the equation
$$(\lambda^2-2)  \varphi(s) - {2\over \sin(\pi s) } \varphi(s) =0,\quad s \in \sigma.\eqno(4.13)  $$
Taking the inverse Mellin transform (1.5) of both sides of the latter equality, we employ the generalized Parseval equality (1.6) and relation (8.4.2.5) in \cite{prud},  Vol. 3 to find
$$(\lambda^2-2) \psi(x)-  {2\over \pi} \int_0^\infty {\psi(t)\over x+t} dt=0,\    x > 0.$$
Thus $f(x)$ by formula (4.2) is a solution of integral equation (4.1) for all $\varphi(s)$ under condition $\varphi (s)= \varphi (1-s)$ such that its inverse Mellin transform satisfies integral equation (4.5). The absolute convergence of the corresponding integral follows from the Schwarz inequality.

\end{proof}

{\bf Corollary 1}. {\it Let $\lambda \in (-\sqrt 2, \sqrt 2)$.  Then the only trivial solution satisfies integral equation $(4.1)$.

\begin{proof} In fact, since $\lambda^2- 2-   2/ \sin(\pi s)  <  0, \ s \in \sigma $, we have from (4.13) $\varphi (s)\equiv 0$ on $\sigma$. Therefore from (4.3) it follows $f^*(s) \equiv 0$ and the inverse Mellin transform implies  $f=0$, i.e.  the solution of (4.1) is trivial. 
\end{proof} 

\noindent {{\bf Acknowledgments}}\\
The present investigation was supported, in part,  by the "Centro de
Matem{\'a}tica" of the University of Porto.\\

\end{document}